\documentclass[amstex,16pt, amssymb]{article}

\usepackage{mathtext}
\usepackage[cp1251]{inputenc}
\usepackage[T2A]{fontenc}
\usepackage[dvips]{graphicx}
\usepackage{amsmath}
\usepackage{amssymb}
\usepackage{amsxtra}
\usepackage{latexsym}
\usepackage{ifthen}
\usepackage{amsthm}

\numberwithin{equation}{section}

\begin{document}

\title{The Dirichlet problem for semi-linear
equations}

\author{Vladimir Gutlyanski\u\i{}, Olga Nesmelova, Vladimir Ryazanov}

\date{}

\maketitle

\begin{abstract}
We study the Dirichlet problem for the semi--linear partial
differential equations ${\rm div}\,(A\nabla u)=f(u)$ in simply
connected domains $D$ of the complex plane $\mathbb C$ with
continuous boundary data. We prove the existence of the weak
solutions $u$ in the class $C\cap W^{1,2}_{\rm loc}(D)$ if a Jordan
domain $D$ satisfies the quasihyperbolic boundary condition by
Gehring--Martio. An example of such a domain that fails to satisfy
the standard (A)--condition by Ladyzhenskaya--Ural'tseva and the
known outer cone condition is given. We also extend our results to
simply connected non-Jordan domains formulated in terms of the prime
ends by Caratheodory.

Our approach is based on the theory of the logarithmic potential,
singular integrals, the Leray--Schauder technique  and a
factorization theorem in \cite{GNR2017}. This theo\-rem allows us to
represent $u$ in the form $u=U\circ\omega,$ where $\omega(z)$ stands
for a quasiconformal mapping of $D$ onto the unit disk ${\mathbb
D}$, ge\-ne\-ra\-ted by the measurable matrix function $A(z),$ and
$U$ is a  solution of the corresponding quasilinear Poisson equation
in the unit disk ${\mathbb D}$. In the end, we give some
applications of these results to various processes of diffusion and
absorption in anisotropic and inhomogeneous media.\footnote{This
work was partially supported by grants of Ministry of Education and
Science of Ukraine, project number is 0119U100421.}
\end{abstract}



\bigskip
\bigskip

{\bf Dedicated to the memory of Professor Bogdan Bojarski for his
great contribution to the theory of quasiformal mappings.}
\bigskip

\normalsize \baselineskip=18.5pt

\vskip 1cm


\section{Introduction}

\medskip

Given a domain $D$ in ${\Bbb C},$ denote by $M^{2\times 2}_K(D)$ the
class of all $2\times 2$ symmetric matrix function
$A(z)=\{a_{jk}(z)\}$ with measurable entries and ${\rm
det}\,A(z)=1,$ satisfying the uniform ellipticity condition
\begin{equation} \label{eqM} {1\over K}\ |\xi|^2\ \leq\ \langle\,
A(z)\, \xi,\,\xi\, \rangle\ \leq\ K\, |\xi|^2 \,\,\, \,\,\, \,\,\,
\mbox{a.e. in}\,\, D \end{equation} for every $\xi\in{\Bbb C}$ where
$1\leq K<\infty$. Further we study the semilinear equations
\begin{equation} \label{eqQUASI} {\rm div\,}[\, A(z)\, \nabla u(z)\, ]\ =\
f(u(z)), \,\,\,z\in D
\end{equation} with continuous functions $f:\mathbb R\to\mathbb R$
either bounded or     $f(t)/t\to 0$ as $t\to\infty$ which describe
many physical phenomena in anisotropic inhomogeneous media.

\medskip

The equations (\ref{eqQUASI}) are closely relevant to the so--called
Beltrami equations. Let $\mu: D\to\mathbb C$ be a measurable
function with $|\mu(z)|<1$ a.e. The equation
\begin{equation}\label{1}
\omega_{\bar{z}}=\mu(z)\cdot \omega_z\
\end{equation}
where $\omega_{\bar z}=(\omega_x+i\omega_y)/2$,
$\omega_{z}=(\omega_x-i\omega_y)/2$, $z=x+iy$, $\omega_x$ and
$\omega_y$ are partial derivatives of the function $\omega$ in $x$
and $y$, respectively, is said to be a {\bf Beltrami equation}. The
equation~\eqref{1} is said to be {\bf nondegenerate} if
$||\mu||_{\infty}<1$. The homeomorphic solutions of nondegenerate
Beltrami's equations~\eqref{1}  with all the first generalized
derivatives by Sobolev are called {\bf quasiconformal mappings}, see
e.g. \cite{Ahlfors:book} and \cite{LV:book}.

We say that a quasiconformal mapping $\omega$ satisfying (\ref{1})
is {\bf agreed with} $A\in M^{2\times 2}_K(D)$ if
\begin{equation}\label{mu} \mu(z)\ =\ \frac{a_{22}(z)-a_{11}(z)-2i
a_{12}(z)}{{\rm det}\,\,(I+A(z))}
\end{equation} where $I$ is the unit $2\times 2$ matrix. Condition (\ref{eqM})
is now written as \begin{equation}\label{ellipticity1} |\mu(z)|\
\leq\ {K-1\over K+1} \,\,\,\,\,\mbox{a.e. in $D$}\ . \end{equation}
Vice versa, given a measurable function $\mu: D\to\mathbb C$,
satisfying (\ref{ellipticity1}), one can invert the algebraic system
(\ref{mu}) to obtain the matrix function $A\in M^{2\times 2}_K(D)$:
\begin{equation}\label{matrix}
A(z)\ =\ \left(\begin{array}{ccc} {|1-\mu|^2\over 1-|\mu|^2}  & {-2{\rm Im}\,\mu\over 1-|\mu|^2} \\
                            {-2{\rm Im}\,\mu\over 1-|\mu|^2}          & {|1+\mu|^2\over 1-|\mu|^2}
                              \end{array}\right).
 \end{equation}
By the existence theorem for (\ref{1}), see e.g. Theorem V.B.3 in
\cite{Ahlfors:book} and Theorem V.1.3 in \cite{LV:book}, any $A\in
M^{2\times 2}_K(D)$ generates a quasiconfomal mapping
$\omega:D\to\mathbb D$.

We also would like to pay attention to a strong interaction between
linear and non-linear elliptic systems in the plane and
quasiconformal mappings. The most general first order linear
homogeneous elliptic system with real coefficients can be written in
the form $f_{\bar z}+\mu(z)f_z +\nu(z) \overline{f_z}=0,$ with
measurable coefficients $\mu$ and $\nu$ such that $|\mu|+|\nu|\leq
(K-1)/(K+1)<1.$ This equation is a particular case of a non-linear
first order system $f_{\bar z}=H(z,f_z)$ where $H:G\times{\Bbb
C}\to{\Bbb C}$ is Lipschitz in the second variable,
$$|H(z,w_1)-H(z,w_2)|\leq {K-1\over K+1}|w_1-w_2|,\,\,\,\, H(z,0)\equiv 0.$$
The principal feature of the above equation is that the difference
of two solutions  need not solve the same equation but each solution
can be represented as {\it a composition of a quasiconformal
homeomorphism and an analytic function.} Thus quasiconformal
mappings become the central tool for the study of these non-linear
systems. A rather comprehensive treatment of the present state of
the theory is given in the excellent book of Astala, Iwaniec and
Martin \cite{AIM-book}. This book contains also an exhaustive
bibliography on the topic.
 In particular, the following fundamental Harmonic
Factorization Theorem for the uniformly elliptic divergence
equations \begin{equation}\label{de} {\rm div\,}A(z,\nabla
u)=0,\,\,\,z\in\Omega, \end{equation} holds, see \cite{AIM-book},
Theorem 16.2.1: Every solution $u\in W^{1,2}_{\rm loc}(\Omega)$
 of the  equation (\ref{de})
can be expressed as {\it the composition  $u(z) = h(f(z))$ of a
quasiconformal homeomorphism $f : \Omega\to G$ and a suitable
harmonic function $h$ on $G.$}

The main goal of this paper is  to point out another application of
quasiconformal mappings to the study of some  semi-linear partial
differential equations, linear part of which contains the elliptic
operator in the divergence form ${\rm div}\,[A(z)\nabla u(z)].$

A fundamental role in the study of the posed problem will play
  Theorem 4.1 in \cite{GNR2017}, that can be considered as a suitable counterpart to the mentioned
above Factorization theorem: a function $u:D\to\mathbb R$ is a weak
solution of (\ref{eqQUASI}) in the class $C\cap W^{1,2}_{\rm
loc}(D)$ if and only if $u=U\circ\omega$ where $\omega :D\to\mathbb
D$ is a quasiconformal mapping agreed with $A$ and $U$ is a weak
solution in the class $C\cap W^{1,2}_{\rm loc}(\mathbb D)$ of the
quasilinear Poisson equation
\begin{equation}\label{QUASII}
\triangle\, U(w)\ =\ J(w)\cdot f(U(w))\ ,\ \ \ \ w\in\mathbb D\ ,
\end{equation}
$J$ denotes the Jacobian of the inverse quasiconformal mapping
$\omega^{-1}:\mathbb D\to D$.

Note that the mapping $\omega^* :=\omega^{-1}$ is extended to a
quasiconformal mapping of $\mathbb C$ onto itself if $\partial D$ is
the so--called quasicircle, see e.g. Theorem II.8.3 in
\cite{LV:book}. By one of the main Bojarski results, see \cite{Bo},
the generalized derivatives of quasiconformal mappings in the plane
are locally integrable with some power $q>2$. Note also that its
Jacobian $J(w)=|\omega^*_w|^2-|\omega^*_{\bar{w}}|^2$, see e.g.
I.A(9) in \cite{Ahlfors:book}. Consequently, in this case $J\in
L^p(\mathbb D)$ for some $p>1$.

In this connection, recall that the image of the unit disk $\mathbb
D$ under a quasiconformal mapping of $\mathbb C$ onto itself is
called a {\bf quasidisk} and its boundary is called a {\bf
quasicircle} or a {\bf quasiconformal curve}. Recall also that a
{\bf Jordan curve} is a continuous one-to-one image of the unit
circle in $\mathbb C$. As known, such a smooth ($C^1$) or Lipschitz
curve is a quasiconformal curve and, at the same time,
quasiconformal curves can be even locally non--rectifiable as it
follows from the well-known Van Koch snowflake example, see e.g. the
point II.8.10 in \cite{LV:book}. The recent book \cite{GH} contains
a comprehensive discussion and numerous characterizations of
quasidisks, see also \cite{Ahlfors:book}, \cite{Ge} and
\cite{LV:book}.

By Theorem 4.7 in \cite{AK}, cf. also Theorem 1 and Corollary in
\cite{BP}, the Jacobian of a qua\-si\-con\-for\-mal homeomorphism
$\omega^*: \mathbb D\to D$ is in $L^p(\mathbb D)$, $p>1$, iff $D$
satisfies the {\bf quasihyperbolic boundary condition}, i.e.
\begin{equation} \label{eqHYPERB}
k_D(z,z_0)\ \le\ a\cdot\ln \frac{d(z_0,\partial D)}{d(z ,\partial
D)}\, +\, b\ \ \ \ \ \ \ \ \ \forall\ z\in D
\end{equation}
for some constants $a$ and $b$ and a fixed point $z_0\in D$ where
$k_D(z,z_0)$ is the {\bf quasi\-hyperbolic distance} between the
points $z$ and $z_0$ in the domain $D$,
\begin{equation} \label{eqHYPERD}
k_D(z,z_0)\ :=\ \inf\limits_{\gamma} \int\limits_{\gamma}
\frac{ds}{d(\zeta,\partial D)}\ .
\end{equation}
Here $d(\zeta,\partial D)$ denotes the Euclidean distance from a
point $\zeta\in D$ to the boun\-da\-ry of $D$ and the infimum is
taken over all rectifiable curves $\gamma$ joining the points $z$
and $z_0$ in $D$.

The notion of domains with the quasihyperbolic boundary condition
was introduced in \cite{GM2} but, before it, was first applied in
\cite{BP}. Note that such domains can be not satisfying the
(A)--condition as well as the outer cone condition and not Jordan at
all, see Sections 4 and 5. Note also that, generally speaking not
rectifiable, quasidisks and, in particular, smooth and Lipschitz
domains satisfy the quasihyperbolic boundary condition.

In Section 2 we give the necessary backgrounds for the Poisson
equation $\triangle\, u(z)=g(z)$ due to the theory of the Newtonian
potential and the theory of singular integrals in $\mathbb C$. First
of all, correspondingly to the key fact of the potential theory,
Proposition 1, the Newtonian potential
\begin{equation}
\label{NEWTON} N_{g}(z)\  :=\ \frac{1}{2\pi}\int\limits_{\mathbb C}
\ln|z-w|\, g(w)\ d\, m(w) \end{equation} of arbitrary integrable
densities $g$ of charge with compact support satisfies the Poisson
equation in a distributional sense, see Corollary 1. Moreover, $N_g$
is continuous for $g\in L^p(\mathbb C)$ and, furthermore, the
Newtonian operator $N:L^p(\mathbb C)\to C(\mathbb C)$ is completely
con\-ti\-nu\-ous for $p>1$, Theorem 1. An examp\-le in Proposition 2
shows that $N_g$ for $g\in L^1(\mathbb C)$ can be not continuous and
even not in $L^{\infty}_{\rm loc}(\mathbb C)$. Theorem 2 says on
additional properties of regularity of $N_g$ depending on a degree
of in\-te\-gra\-bi\-li\-ty of $g$. Finally, resulting Corollary 2
states the existence, representation and regularity of solutions to
the Dirichlet problem for the Poisson equation with continuous
boundary data.

Section 3 contains one of  the main results Theorem 3 on the
existence of regular solutions of the Dirichlet problem to the
quasilinear Poisson equation
\begin{equation}\label{Poisson1} \triangle u(z)\ =\ h(z)\cdot
f(u(z))
\end{equation} in the unit disk $\mathbb D$ for
arbitrary continuous boundary data. In general, we assume that the
function $h:\mathbb D\to\mathbb R$ is in the class $L^p(\mathbb D)$,
$p>1,$ and the continuous function $f:\mathbb R\to\mathbb R$ is
either bounded or with non--decreasing $|f\,|$ of $\ |t|$ and
\begin{equation}\label{lim}
\lim_{t\to \infty}{f(t)\over t}\ =\ 0\ ,
\end{equation}
without any assumptions on the sign and zeros of the right hand side
in (\ref{Poisson1}). The degree of regularity of solutions depends
first of all on a degree of in\-te\-gra\-bi\-li\-ty of the
multiplier $h$. The proof of Theorem 3 is realized by the
Leray--Schauder approach. This result is extended to arbitrary
smooth ($C^1$) domains in Corollary 3. The proof of the latter is
obtained from Theorem 3 through the fundamental results of
Caratheodory--Osgood--Taylor and Warschawski on the boundary
behavior of conformal mappings between Jordan domains in $\mathbb
C$.

Section 4 includes the main result of the present paper Theorem 4 on
the existence of regular weak solutions of the Dirichlet problem for
semi-linear equations (\ref{eqQUASI}) in Jordan domains $D$ in
$\mathbb C$ satisfying the quasihyperbolic boun\-da\-ry condition
with arbitrary con\-ti\-nu\-ous boundary data $\varphi :\partial
D\to\mathbb R$. Theorem 4 states the existence of a weak solution
$u:{D}\to\mathbb R$ of the equation (\ref{eqQUASI}) in the class
$C\cap W^{1,2}_{\rm loc}(\Omega)$ which is locally H\"older
continuous in $D$ and continuous in $\overline{D}$ with
$u|_{\partial D}=\varphi$. If in addition $\varphi$ is H\"older
continuous, then $u$ is H\"older continuous in $\overline{ D}$.
Moreover, $u=U\circ\omega$ where $U$ is a weak solution of the
quasilinear Poisson equation (\ref{QUASII}) and $\omega :D\to\mathbb
D$ is a quasiconformal mapping agreed with $A$. In Lemma 2 we show
that there exist Jordan domains $D$ in $\mathbb C$ with the
quasihyperbolic boundary condition but without the standard
(A)--condition and, consequently, without the known outer cone
condition.


Section 5 contains Theorem 5 which is the extension of Theorem 4 to
arbitrary bounded simple connected (not Jordan !) domains $D$ in
$\mathbb C$ with the quasihyperbolic boun\-da\-ry condition
formulated in terms of the prime ends by Caratheodory. Finally, in
Section 6 we give applications of the obtained results to various
kind of absorption, reaction-diffusion problems, equations of a
heated plasma and the combustion in anisotropic and inhomogeneous
media, Theorem 6 and 7, see also Corollaries 6--8.


Later on, as usual $\ C^{\infty}_0(\mathbb C)$ denotes the class of
all infinitely differentiable functions $\psi :\ \mathbb C\to
\mathbb R\ $ with compact support,
$\triangle=\frac{\partial^2}{\partial
x^2}+\frac{\partial^2}{\partial y^2}$ is the Laplace differential
operator and $dm(z)=dxdy$, $z=x+iy$, corresponds to the Lebesgue
measure and $\mathbb D:=\{ z\in\mathbb C : |z|<1\}$ is the unit disk
in the complex plane $\mathbb C$.

\section{Potentials and the Poisson equation}

Let us recall the {\bf Green function} and the {\bf Poisson kernel}
in the unit disk $\mathbb D$
\begin{equation}
\label{GP} G(z,w) := \ln \left| \frac{1 - z\bar{w}}{z - w} \right|\
,\  z\ne w,\ \ \ \ \hbox{and}\ \ \ P(z, e^{it}) := \frac{1 - |z|^2}
{|1 - ze^{-it}|^2}\end{equation} and, for real-valued functions
$\varphi\in C(\partial\mathbb D)$ and $g\in C(\overline{\mathbb
D})$, apply the notations
\begin{equation} \label{operators} {\cal{P}}_{\varphi} (z)\ =\
\frac{1}{2\pi}\ \int\limits_0^{2\pi} P(z, e^{it})\,
\varphi(e^{-it})\, dt\ ,\ \ \ {\cal{G}}_g(z)\ =\
\int\limits_{\mathbb D} G(z,w)\, g(w)\, dm(w)\ .
\end{equation}
As known, see e.g. \cite{Ho}, p. 118-120, then a solution to the
{\bf Poisson equation}
\begin{equation}
\label{Poisson} \triangle f(z)\  =\ g(z)\ ,\ \ \ z\in\mathbb D\
,\end{equation} satisfying the boundary condition
$f|_{\partial\mathbb D} = \varphi$, is given by the formula
\begin{equation}
\label{solution}  f(z)\ =\ {\cal{P}}_{\varphi} (z) - {\cal{G}}_g(z)\
.\end{equation}

In this section, we give the representation of solutions of the
Poisson equation in the form of the Newtonian (normalized
antilogarithmic) potential that is more convenient for our research
and, on this basis, we prove the existence and representation
theorem for solutions of the Dirichlet problem to the Poisson
equation under the corresponding conditions of integrability of
sources $g$.

Correspondingly to 3.1.1 in \cite{Ra}, given a finite Borel measure
$\nu$ on $\mathbb C$ with compact support, its {\bf potential} is
the function $p_{\nu}:\mathbb C\to[-\infty,\infty)$ defined by
\begin{equation}
\label{eqPOTENTIAL} p_{\nu}(z)\  =\ \int\limits_{\mathbb C}
\ln|z-w|\, d\nu(w)\ .\end{equation}

{\bf Remark 1.} Note that the function $p_{\nu}$ is subharmonic by
Theorem 3.1.2 and, consequently, it is locally integrable on
$\mathbb C$ by Theorem 2.5.1 in \cite{Ra}. Moreover,  $p_{\nu}$ is
harmonic outside the support of $\nu$.


This definition can be extended to finite {\bf charges}  $\nu$ with
compact support (named also {\bf signed measures}), i.e., to real
valued sigma-additive functions on Borel sets in $\mathbb C$,
because of $\nu = \nu^+-\nu^-$ where $\nu^+$ and $\nu^-$ are Borel
measures by the well--known Jordan decomposition, see e.g. Theorem
0.1 in \cite{La}.

\bigskip

The key fact is the following statement, see e.g. Theorem 3.7.4 in
\cite{Ra}.

\medskip

{\bf Proposition 1.} {\it Let $\nu$ be a finite charge with compact
support in $\mathbb C$. Then
\begin{equation}
\label{eqLAPLACE} \triangle p_{\nu}\  =\ 2\pi\cdot\nu\end{equation}
in the distributional sense, i.e.,
\begin{equation}
\label{eqDISTRIBUTION} \int\limits_{\mathbb C}
p_{\nu}(z)\,\triangle\psi(z)\, d\,m(z)\ =\ 2\pi\int\limits_{\mathbb
C} \psi(z)\, d\,\nu(z)\ \ \ \ \ \ \ \ \  \forall\ \psi\in
C^{\infty}_0(\mathbb C)\ .\end{equation}}

If $g:\mathbb C\to\mathbb R$ is an integrable function with compact
support and
\begin{equation} \label{eqINTEGRABLE} \nu(B)\ :=\ \int\limits_Bg(z)\ d\,m(z)\end{equation}
for every Borel set $B$ in $\mathbb C$, then $g$ is called a {\bf
density of charge} $\nu$ and the {\bf Newtonian potential} of $g$ is
the function
\begin{equation} \label{eqIPOTENTIAL} N_{g}(z)\  :=\
\frac{1}{2\pi}\int\limits_{\mathbb C} \ln|z-w|\, g(w)\ d\, m(w)\
.\end{equation}

{\bf Corollary 1.} {\it If $g:\mathbb C\to\mathbb R$ is an
integrable function with compact support, then
\begin{equation}
\label{eqILAPLACE} \triangle N_{g}\  =\  g\end{equation} in the
distributional sense, i.e.,
\begin{equation}
\label{eqIDISTRIBUTION} \int\limits_{\mathbb C}
N_{g}(z)\,\triangle\psi(z)\, d\,m(z)\ =\ \int\limits_{\mathbb C}
\psi(z)\, g(z)\, d\,m(z)\ \ \ \ \ \ \ \ \  \forall\ \psi\in
C^{\infty}_0(\mathbb C)\ .\end{equation}}

The {\bf shift of a set $E\subset\mathbb C$ by a complex vector
$\Delta z\in\mathbb C$} is the set
$$
E+\Delta z\ :=\ \{\ \xi\in\mathbb C:\ \xi\ =\ z\ +\ \Delta z\ ,\
z\in E\ \}\ .
$$
The next statement on continuity in the mean of functions $\psi
:\mathbb C\to\mathbb R$  in $L^q(\mathbb C)$, $q\in[1,\infty)$, with
respect to shifts is useful for the study of the Newtonian
potential, see e.g. Theorem 1.4.3 in \cite{So}, cf. also Theorem
III(11.2) in \cite{Sa}. The one-dimensional analog of the statement
can be found also in \cite{Ru}, Theorem 9.5.

\medskip

{\bf Lemma 1.} {\it Let $\psi\in L^q(\mathbb C)$, $q\in[1,\infty)$,
have a compact support. Then \begin{equation} \label{eqSHIFTS}
\lim\limits_{\Delta z\to 0}\int\limits_{\mathbb C}|\psi(z+\Delta
z)-\psi(z)|^q\ d\, m(z)\ =\ 0\ .
\end{equation}}

Here we  give an alternative  proof of this important statement that
may be of independent interest. The proof is based on arguments by
contradiction and the absolute continuity of indefinite integrals.

{\bf Proof.} Let us assume that there is a sequence $\Delta
z_n\in\mathbb C$, $n=1,2,\ldots$, such that $\Delta z_n\to 0$ as
$n\to\infty$ and, for some $\delta >0$ and $\psi_n(z):=\psi(z+\Delta
z_n)$, $n=1,2,\ldots$,
\begin{equation}
\label{eqCONTRA} I_n\ :=\ \left[\int\limits_{\mathbb
C}|\psi_n(z)-\psi(z)|^q\ d\, m(z)\right]^{\frac{1}{q}}\ \ge\ \delta\
\ \ \ \ \ \forall\ n=1,2,\ldots\ .
\end{equation}

Denote by $K$ the closed disk in $\mathbb C$ centered at $0$ with
the minimal radius $R$ that contains the support of $\psi$. By the
Luzin theorem, see e.g. Theorem 2.3.5 in \cite{Fe}, for every
prescribed $\varepsilon>0$, there is a compact set $C\subset K$ such
that $g|_C$ is continuous and $m(K\setminus C)<\varepsilon$. With no
loss of generality, we may assume that $C\subset K_*$ where $K_*$ is
a closed disk in $\mathbb C$ centered at $0$ with a radius
$r\in(0,R)$ and, moreover, that $C_n\subset K$, where $C_n:=C-\Delta
z_n$, for all $n=1,2,\ldots $. Note that $m(C_n)=m(C)$ and then
$m(K\setminus C_n)<\varepsilon$ and, consequently, $m(K\setminus
C^*_n)<2\varepsilon$, where $C_n^*:=C\cap C_n$, because $K\setminus
C^*_n=(K\setminus C_n)\cup(K\setminus C)$.

Next, setting $K_n=K-\Delta z_n$, we see that $K\cup
K_n=C_n^*\cup(K\setminus C_n^*)\cup(K_n\setminus C_n^*)$ and that
$K_n\setminus C_n^*+\Delta z_n=K\setminus C_n^*$. Hence by the
triangle inequality for the norm in $L^p$ the following estimate
holds $$ I_n\ \le\ 4\ \cdot \left[\int\limits_{K\setminus
C^*_n}|\psi(z)|^q\ d\, m(z)\right]^{\frac{1}{q}}\ +\
\left[\int\limits_{C^*_n}|\psi_n(z)-\psi(z)|^q\ d\,
m(z)\right]^{\frac{1}{q}} \ \ \  \forall\ n=1,2,\ldots $$ By
construction the both terms from the right hand side can be made to
be arbitrarily small, the first one for small enough $\varepsilon$
because of absolute continuity of indefinite integrals and the
second one for all large enough $n$ after the choice of the set $C$.
Thus, the assumption (\ref{eqCONTRA}) is disproved. $\Box$

\medskip

{\bf Theorem 1.} {\it Let $g:\mathbb C\to\mathbb R$ be in
$L^p(\mathbb C)$, $p>1$, with compact support. Then $N_g$ is
continuous. A collection $\{N_g\}$ is equicontinuous if the
collection $\{ g\}$ is bounded by the norm in $L^p(\mathbb C)$ with
supports in a fixed disk $K$.  Moreover, under these hypothesis, on
each compact set $S$ in $\mathbb C$
\begin{equation} \label{eqEST}
\| N_g\|_C\ \le\ M\cdot \| g\|_p\ ,
\end{equation} where $M$ is a constant depending in general on $S$ but not on $g$. }

The corresponding statement on the continuity of integrals of
potential type in $\mathbb R^n$, $n\ge 3$, can be found in
\cite{So}, Theorem 1.6.1.

\medskip

{\bf Proof.} By the H\"older inequality with
$\frac{1}{q}+\frac{1}{p}=1$ we have that
$$
|N_{g}(z)-N_{g}(\zeta)|\ \le\ \frac{\| g\|_p}{2\pi}\cdot \left[\
\int\limits_{K} |\, \ln |z-w|-\ln |\zeta-w|\, |^q\ d m(w)\,
\right]^{\frac{1}{q}}\ =
$$
$$
=\ \frac{\| g\|_p}{2\pi}\cdot \left[\ \int\limits_{\mathbb C} |\,
\psi_{\zeta}(\xi+\Delta z)-\psi_{\zeta}(\xi)\, |^q\ d m(\xi)\,
\right]^{\frac{1}{q}}$$ where $\xi =\zeta -w$, $\Delta z=z-\zeta$,
$\psi_{\zeta}(\xi):=\chi_{K+\zeta}(\xi)\ln |\xi|$. Thus, the first
two conclusions follow by Lemma 1 because of the function $\ln
|\xi|$ belongs to the class $L^q_{\rm loc}(\mathbb C)$ for all $q\in
[1,\infty)$.

The third conclusion similarly follows  through the direct estimate
$$
|N_{g}(\zeta)|\ \le\ \frac{\| g\|_p}{2\pi}\cdot \left[\
\int\limits_{K} |\, \ln |\zeta-w|\, |^q\ d m(w)\,
\right]^{\frac{1}{q}}  =\ \ \frac{\| g\|_p}{2\pi}\cdot \left[\
\int\limits_{\mathbb C} |\, \psi_{\zeta}(\xi)\, |^q\ d m(\xi)\,
\right]^{\frac{1}{q}} $$ because the latter integral is continuous
in $\zeta\in\mathbb C$. Indeed, by the triangle ine\-qua\-li\-ty for
the norm in  $L^q(\mathbb C)$ we see that
$$ \left| \| \psi_{\zeta}\|_q\ -\ \| \psi_{\zeta_*}\|_q  \right|\ \le\ \| \psi_{\zeta}-\psi_{\zeta_*}\|_q\ =\ \left\{\
\int\limits_{\Delta} |\ln |\xi||^q\ d m(\xi)\,
\right\}^{\frac{1}{q}}
$$
where $\Delta$ denotes the symmetric difference of the disks
$K+\zeta$ and $K+\zeta_*$. Thus, the statement follows from the
absolute continuity of the indefinite integral. $\Box$

\medskip

{\bf Proposition 2.} {\it There exist functions $g\in L^1(\mathbb
C)$ with compact support whose potentials $N_g$ are not continuous,
furthermore, $N_g\notin L^{\infty}_{\rm loc}$.}

\bigskip

{\bf Proof.} Indeed, let us consider the function
$$\omega(t):=\frac{1}{t^2(1-\ln t)^{\alpha}}\ ,\ t\in(0,1]\ ,\ \alpha\in(1,2)\ ,\ \ \ \ \ \ \omega(0):=\infty$$
and, correspondingly,
$$g(z):=\omega(|z|)\ ,\ z\in\overline{\mathbb D}\ ,\ \ \ \ \ \ \  g(z):= 0\ ,\
z\in\mathbb C\setminus\overline{\mathbb D}\ .$$

\bigskip

Then, setting $\Omega(t)=t\cdot\omega(t)$, we see, firstly, that
$$
\int\limits_{\overline{\mathbb D}}|g(w)|\ d\, m(w)\ =\
2\pi\lim_{\rho\to +0} \int\limits_{\rho}^{1}\Omega(t)\ d\, t\ =\
2\pi\lim_{\rho\to +0}\int\limits_{\rho}^{1}\frac{d \ln t}{(1-\ln
t)^{\alpha}}\ =\ \frac{2\pi}{\alpha -1}
$$
and, secondly, that the Newtonian potential $N_g$ at the origin is
equal to
$$
\lim_{\rho\to +0}\int\limits_{\rho}^{1}\Omega(t)\ln t\, d\, t =
\lim_{\rho\to +0}\left\{\left[\ln \frac{1}{t}
\int\limits_{t}^{1}\Omega(\tau)\, d\, \tau \right]_{\rho}^{1} +
\int\limits_{\rho}^{1}\left(\frac{1}{t}
\int\limits_{t}^{1}\Omega(\tau)\, d\, \tau\right) d\, t\right\} =
$$
$$
=\ \frac{1}{\alpha -1}\cdot\lim_{\rho\to +0}\left(\left[\frac{\ln
t}{(1-\ln t)^{\alpha - 1}}\right]_{\rho}^{1}\ -\
\int\limits_{\rho}^{1}\frac{d\, t}{t(1-\ln t)^{\alpha - 1}}\right)\
=
$$
$$=\ \frac{1}{\alpha
-1}\cdot\lim_{\rho\to +0}\left[\frac{1}{(1-\ln t)^{\alpha -1}}\
+\frac{\alpha -1}{2-\alpha}\cdot (1-\ln t)^{2-\alpha}
\right]_{\rho}^{1}\ =\ -\infty\ .\ \ \ \Box
$$

\bigskip

Next, recall the definition of the formal complex derivatives:
$$
\frac{\partial}{\partial z}\ :=\ \frac{1}{2}\left\{
\frac{\partial}{\partial x}\ -\ i\cdot\frac{\partial}{\partial
y}\right\}\ , \ \ \ \frac{\partial}{\partial\overline z}\ :=\
\frac{1}{2}\left\{ \frac{\partial}{\partial x}\ +\
i\cdot\frac{\partial}{\partial y}\right\}\ ,\ \ \ z=x+iy\ .
$$

\bigskip

The elementary algebraic calculations show their relation to the
Laplacian
$$
\triangle\ :=\ \frac{\partial^2}{\partial x^2}\ +\
\frac{\partial^2}{\partial y^2}\ =\ 4\cdot\frac{\partial
^2}{\partial z\partial\overline z}\ =\ 4\cdot\frac{\partial
^2}{\partial\overline z\partial z}\ .
$$

\bigskip


Further we apply the theory of the well-known integral operators
$$
Tg(z)\ :=\ \frac{1}{\pi}\int\limits_{\mathbb C} g(w)\ \frac{d\,
m(w)}{z-w}\ ,\ \ \ \ \overline Tg(z)\ :=\
\frac{1}{\pi}\int\limits_{\mathbb C} g(w)\ \frac{d\, m(w)}{\overline
z-\overline w}
$$
defined for integrable functions with a compact support $K$ and
studied in detail. Recall the known results on them in Chapter 1 of
\cite{Ve}, confining the case $K=\overline{\mathbb D}$, that are
relevant to the proof of Theorem 2.

First of all,  if $g\in L^1(\mathbb C)$, then by Theorem 1.13 the
integrals $Tg$ and $\overline Tg$ exist a.e. in $\mathbb C$ and
belong to $L^q_{\rm loc}(\mathbb C)$ for all $q\in[1,2)$ and by
Theorem 1.14 they have generalized derivatives by Sobolev
$(Tg)_{\overline z}=g=(\overline Tg)_z$. Furthermore, if  $g\in
L^p(\mathbb C)$, $p>1$, then by Theorem 1.27 and (6.27) $Tg$ and
$\overline Tg$ belong to $L^{q}_{\rm loc}(\mathbb C)$ for some
$q>2$. and, moreover, by Theorems 1.36--1.37 $(Tg)_{z}$ and
$(\overline Tg)_{\overline z}$ also belong to $L^{p}_{\rm
loc}(\mathbb C)$. Finally, if $g\in L^p(\mathbb C)$ for $p>2$, then
by Theorem 1.19  $Tg$ and $\overline Tg$ belong to $C^{\alpha}_{\rm
loc}(\mathbb C)$ with $\alpha = (p-2)/p$.

\bigskip

The following theorem on the Newtonian potentials is important to
obtain solutions of the Dirichlet problem to the Poisson equation of
a higher regularity.

\medskip

{\bf Theorem 2.} {\it Let $g:\mathbb C\to\mathbb R$ have compact
support. If $g\in L^1(\mathbb C)$, then $N_g\in L^r_{\rm loc}$ for
all $r\in[1,\infty)$, $N_g\in W^{1,q}_{\rm loc}$ for all
$q\in[1,2)$. Moreover, there exist generalized derivatives by
Sobolev $\frac{\partial^2 N_{g}}{\partial z\partial\overline z}$ and
$\frac{\partial^2 N_{g}}{\partial\overline z\partial z}$ and
\begin{equation}\label{eqLAP}
4\cdot\frac{\partial^2 N_{g}}{\partial z\partial\overline z}\ =\
\triangle N_g\ =\ 4\cdot\frac{\partial^2 N_{g}}{\partial\overline
z\partial z}\ =\ g\ \ \ \mbox{a.e.} \end{equation} If $g\in
L^p(\mathbb C)$, $p>1$, then $N_g\in W^{2,p}_{\rm loc}$, $\triangle
N_g=g$ a.e. and, moreover, $N_g\in W^{1,q}_{\rm loc}$ for $q>2$,
consequently, $N_g$ is locally H\"older continuous. If $g\in
L^p(\mathbb C)$, $p>2$, then $N_g\in C^{1,\alpha}_{\rm loc}$ where
$\alpha = (p-2)/p$.}

\medskip

{\bf Proof.} Note that $N_g$ is the convolution $\psi *g$, where
$\psi(\zeta)=\ln |\zeta|$, and hence $N_g\in L^r_{\rm loc}$ for all
$r\in[1,\infty)$, see e.g. Corollary 4.5.2 in \cite{Hor}. Moreover,
as well--known $\frac{\partial\psi
*g}{\partial z}=\frac{\partial\psi
}{\partial z}*g$ and $\frac{\partial\psi
*g}{\partial\overline z}=\frac{\partial\psi
}{\partial\overline z}*g$, see e.g. (4.2.5) in \cite{Hor}, and in
addition by elementary calculations
$$
\frac{\partial}{\partial z}\ln |z-w|\ =\
\frac{1}{2}\cdot\frac{1}{z-w}\ ,\ \ \ \
\frac{\partial}{\partial\overline z}\ln |z-w|\ =\
\frac{1}{2}\cdot\frac{1}{\overline z-\overline w}\ .
$$
Consequently,
$$
\frac{\partial N_{g}(z)}{\partial z}\  =\ \frac{1}{4}\cdot Tg(z)\ ,\
\ \ \ \frac{\partial N_{g}(z)}{\partial\overline z}\  =\
\frac{1}{4}\cdot\overline Tg(z)\ .$$

Thus, all the rest conclusions for $g\in L^1(\mathbb C)$ follow by
Theorems 1.13--1.14 in \cite{Ve}. If $g\in L^p(\mathbb C)$, $p>1$,
then $N_g\in W^{1,q}_{\rm loc}$, $q>2$, by Theorem 1.27, (6.27) in
\cite{Ve}, consequently, $N_g$ is locally H\"older continuous, see
e.g. Lemma 1.7 in \cite{BI}, and $N_g\in W^{2,p}_{\rm loc}$ by
Theorems 1.36--1.37 in \cite{Ve}. If $g\in L^p(\mathbb C)$, $p>2$,
then $N_g\in C^{1,\alpha}_{\rm loc}$ with $\alpha = \frac{p-2}{p}$
by Theorem 1.19 in \cite{Ve}. $\Box$

\medskip

{\bf Remark 2.} Note that generally speaking $N_g\notin W^{2,1}_{\rm
loc}$ for the case $g\in L^1(\mathbb C)$, see e.g. example 7.5 in
\cite{GM$^*$}, p.141. Note also that the corresponding Newtonian
potentials $N_g$ in $\mathbb R^n$, $n\ge 3$, also belong to
$W^{2,p}_{\rm loc}$ if $g\in L^p(\mathbb C)$ for $p>1$ with compact
support, see e.g. \cite{GT}, Theorem 9.9.

\medskip

Further we also apply the following particular case of Theorem 1.10
in \cite{Ve}.

\medskip

{\bf Proposition 3.} {\it Let a function $\varphi : \partial\mathbb
D\to\mathbb R$ be H\"older continuous of the order $\alpha\in(0,1)$.
Then the integral of the Cauchy type
$$
\Phi(z)\ :=\ \frac{1}{2\pi i}\int\limits_{\partial\mathbb
D}\frac{\varphi(\zeta)\, d\zeta}{\zeta -z}\ =\
\frac{1}{2\pi}\int\limits_{0}^{2\pi}\varphi(e^{it})\frac{
e^{it}}{e^{it} -z}\, dt\ ,\ \ \ \ z\in\overline{\mathbb D}\ ,
$$
is also H\"older continuous in $\overline{\mathbb D}$ of the same
order with $\Phi|_{\partial\mathbb D}\equiv \varphi$.}

\medskip

By Theorem 2, Proposition 3 and the known Poisson formula, see e.g.
I.D.2 in \cite{Ko}, we come to the following consequence on the
existence, regularity and representation of solutions for the
Dirichlet problem to the Poisson equation in the unit disk $\mathbb
D$ where we assume the charge density $g$ to be extended by zero
outside $\mathbb D$.

\medskip

{\bf Corollary 2.} {\it Let $\varphi : \partial\mathbb D\to\mathbb
R$ be a continuous function and $g: \mathbb D\to\mathbb R$ belong to
the class $L^p(\mathbb D)$, $p>1$. Then the function $U:=N_g - {\cal
P}_{N^*_g} + {\cal P}_{\varphi}$, $N_g^*:=N_g|_{\partial\mathbb D}$,
is continuous in $\overline{\mathbb D}$ with $U|_{\partial\mathbb
D}=\varphi$, belongs to the class $W^{2,p}_{\rm loc}(\mathbb D)$ and
$\triangle U =g$ a.e. in $\mathbb D$. Moreover, $U\in W^{1,q}_{\rm
loc}(\mathbb D)$ for some $q>2$ and $U$ is locally H\"older
continuous. If in addition $\varphi$ is H\"older continuous, then
$U$ is H\"older continuous in $\overline{\mathbb D}$. If $g\in
L^p(\mathbb D)$, $p>2$, then $U\in C^{1,\alpha}_{\rm loc}(\mathbb
D)$ where $\alpha = (p-2)/p$.}

\medskip

{\bf Remark 2.} The H\"older continuity of $U$ for H\"older
continuous $\varphi$ follows from Proposition 3 because of the
Poisson kernel $$P(z,e^{it})\ =\ {\rm Re}\
\frac{e^{it}+z}{e^{it}-z}\ =\ -1\ +\ 2\,{\rm Re}\
\frac{e^{it}}{e^{it}-z}\ .$$ Note also by the way that a generalized
solution of the Dirichlet problem to the Poisson equation in the
class $C(\overline{\mathbb D} )\cap W^{1,2}_{\rm loc}(\mathbb D)$ is
unique at all, see e.g. Theorem 8.30 in \cite{GT}. One can show that
the integral operators generated derivatives of solutions in Theorem
2 and Corollary 2 are completely continuous (it is clear from the
cor\-res\-pon\-ding theorems in \cite{Ve} mentioned under the proof
of Theorem 2), cf. e.g. \cite{Ka} and \cite{KV}. However, for our
goals it is sufficient that the operator $N_g: L^p({\mathbb D})\to
C(\overline{\mathbb D})$ is completely continuous by Theorem 1 for
$p>1$, see the proof of Theorem 3 further.

\bigskip

\section{The case of the quasilinear Poisson equations}

\ \ \ \ \ The case of the quasilinear Poisson equations will be
reduced to the linear Poisson equation from the previous section by
the Leray--Schauder approach. In this connection, recall the main
relevant content of the classical paper \cite{LS}.

First of all, Leray and Schauder define a {\bf completely
continuous} mapping from a metric space $M_1$ into a metric space
$M_2$ as a continuous mapping on $M_1$ which takes bounded subsets
of $M_1$ into {\bf relatively compact} ones of $M_2$, i.e. with
compact closures in $M_2$. When a continuous mapping takes $M_1$
into a relatively compact subset of $M_1$, it is nowadays said to be
{\bf compact} on $M_1$.

Then Leray and Schauder extend as follows the Brouwer degree to
compact perturbations of the identity $I$ in a Banach space $B$,
i.e. a complete normed linear space. Namely, given an open bounded
set $\Omega\subset B$, a compact mapping $F: B\to B$ and $z \notin
\Phi(\partial \Omega)$, $\Phi :=I-F$, the {\bf (Leray–Schauder)
topological degree} $\deg\, [\Phi,\Omega, z]$ of $\Phi$ in $\Omega$
over $z$ is constructed from the Brouwer degree by approximating the
mapping $F$ over $\Omega$ by mappings $F_{\varepsilon}$ with range
in a finite-dimensional subspace $B_{\varepsilon}$ (containing $z$)
of $B$. It is showing that the Brouwer degrees $\deg\,
[\Phi_{\varepsilon} ,\Omega_{\varepsilon}, z]$ of
$\Phi_{\varepsilon}:=I_{\varepsilon} - F_{\varepsilon}$,
$I_{\varepsilon}:=I|_{B_{\varepsilon}}$, in
$\Omega_{\varepsilon}:=\Omega\cap B_{\varepsilon}$ over $z$
stabilize for sufficiently small positive $\varepsilon$ to a common
value defining $\deg\, [\Phi,\Omega, z]$ of $\Phi$ in $\Omega$ over
$z$.

This topological degree “algebraically counts” the number of fixed
points of $F(\cdot)-z$ in $\Omega$ and conserves the basic
properties of the Brouwer degree as ad\-di\-ti\-vi\-ty and homotopy
invariance. Now, let $a$ be an isolated fixed point of $F$. Then the
{\bf local (Leray–Schauder) index} of $a$ is defined by ${\rm ind}\,
[\Phi, a] := \deg [\Phi,B(a, r), 0]$ for small enough r > 0. If
$a=0$, then we say on the {\bf index} of $F$. In particular, if
$F\equiv 0$, correspondingly, $\Phi\equiv I$, then the index of $F$
is equal to $1$.

{\bf Theorem 1 in \cite{LS}} can be formulated in the following way.
{\it Let $B$ be a Banach space, and let $F(\cdot,\tau):B\to B$ be a
family of operators with $\tau\in[0,1]$. Suppose that the following
hypotheses hold:

{\rm {\bf (H1)}} $F(\cdot,\tau)$ is completely continuous on $B$ for
each $\tau\in[0,1]$ and uniformly continuous with respect to the
parameter $\tau\in[0,1]$ on each bounded set in $B$;

{\rm {\bf (H2)}} the operator $F:=F(\cdot,0)$ has finite collection
of fixed points whose total index is not equal to zero;

{\rm {\bf (H3)}} the collection of all fixed points of the operators
$F(\cdot,\tau)$, $\tau\in[0,1]$,  is bounded in $B$.

Then the collection of all fixed points of the family of operators
$F(\cdot,\tau)$ contains a continuum along which $\tau$ takes all
values in  $[0,1]$.}


In the proof of the next theorem the initial operator
$F(\cdot):=F(\cdot,0)\equiv 0$. Hence $F$ has the only one fixed
point (at the origin) and its index is equal to $1$ and, thus,
hypothesis (H2) will be automatically satisfied.


{\bf Theorem 3.} {\it Let $\varphi :\partial\mathbb D\to\mathbb R$
be a continuous function, $h:\mathbb D\to\mathbb R$ be a function in
the class $L^p(\mathbb D)$, $p>1$,  and a function $f:\mathbb
R\to\mathbb R$ be continuous and
\begin{equation} \label{eqAPRIORY}
\lim\limits_{t\to \infty}\ \frac{f(t)}{t}\ =\ 0\ .
\end{equation} Then there is a continuous function $U:\overline{\mathbb D}\to\mathbb R$
with $U|_{\partial\mathbb D}=\varphi$, $U|_{\mathbb D}\in
W^{2,p}_{\rm loc}$ and
\begin{equation} \label{eqQUASILINEAR}
\triangle\, U(z)\ =\ h(z)\cdot f(U(z))\ \ \ \ \ \ \ \mbox{for a.e.}\
z\in\mathbb D\ .
\end{equation} Moreover, $U\in W^{1,q}_{\rm loc}(\mathbb
D)$ for some $q>2$ and $U$ is locally H\"older continuous in
$\mathbb D$. If in addition $\varphi$ is H\"older continuous, then
$U$ is H\"older continuous in $\overline{\mathbb D}$. Furthermore,
if $p>2$, then $U\in C^{1,\alpha}_{\rm loc}(\mathbb D)$ where
$\alpha = (p-2)/p$. In particular, $U\in C^{1,\alpha}_{\rm
loc}(\mathbb D)$ for all $\alpha\in(0,1)$ if $h\in
L^{\infty}(\mathbb D)$.}

\medskip

{\bf Proof.} If $\| h\|_p=0$ or $\| f\|_C= 0$, then the Poisson
integral ${\cal P}_{\varphi}$ gives the desired solution of the
Dirichlet problem for equation (\ref{eqQUASILINEAR}), see e.g. I.D.2
in \cite{Ko}. Hence we may assume further that  $\| h\|_p\neq 0$ and
$\| f\|_C\neq 0$.

Set $f_*(s)=\max\limits_{|t|\le s}|f(t)|$, $s\in\mathbb
R^+:=[0,\infty)$. Then the function $f_*:\mathbb R^+\to\mathbb R^+$
is continuous and nondecreasing and, moreover, $f_*(s)/s\to 0$ as
$s\to\infty$ by (\ref{eqAPRIORY}).

By Theorem 1 and the maximum principle for harmonic functions, we
obtain the family of operators $F(g;\tau): L^{p}(\mathbb D)\to
L^{p}(\mathbb D)$, $\tau\in[0,1]$:
\begin{equation} \label{eqFORMULA}
F(g;\tau)\ :=\ \tau h\cdot f(N_g - {\cal P}_{N^*_g} + {\cal
P}_{\varphi})\ , \ N_g^*:=N_g|_{\partial\mathbb D}\ ,\ \ \ \ \ \ \
\forall\ \tau\in[0,1]
\end{equation}
which satisfies all groups of hypothesis (H1)-(H3) of Theorem 1 in
\cite{LS}. Indeed:

(H1). First of all, $F(g;\tau)\in L^{p}(\mathbb D)$ for all
$\tau\in[0,1]$ and  $g\in L^{p}(\mathbb D)$ because by Theorem 1
$f(N_g - {\cal P}_{N^*_g} + {\cal P}_{\varphi})$ is a continuous
function and, moreover, by (\ref{eqEST}) $$\| F(g;\tau)\|_p\ \le\ \|
h\|_p\ f_*\left(\, 2M\, \| g\|_p + \| \varphi\|_C\right) \ <\
\infty\ \ \ \ \ \forall\ \tau\in[0,1]\ .$$ Thus, by Theorem 1 in
combination with the Arzela--Ascoli theorem, see e.g. Theorem IV.6.7
in \cite{DS}, the operators $F(g;\tau)$ are completely continuous
for each $\tau\in[0,1]$ and even uniformly continuous with respect
to the parameter $\tau\in[0,1]$.

(H2). The index of the operator $F(g;0)$ is obviously equal to $1$.

(H3). By inequality (\ref{eqEST}) and the maximum principle for
harmonic functions, we have the estimate for solutions $g\in L^p$ of
the equations $g=F(g;\tau)$:
$$
\| g\|_p\ \le\ \| h\|_p\ f_*\left(\, 2M\, \| g\|_p + \|
\varphi\|_C\right) \ \le\ \| h\|_p\ f_*(\, 3M\, \| g\|_p) \ \ \ \ \
\ \ \
$$
whenever $\|g\|_p\ge \| \varphi\|_C/M$, i.e. then it should be
\begin{equation} \label{eqEST3}
\frac{f_*(\, 3M\, \| g\|_p)}{3M\, \| g\|_p}\ \ge\ \frac{1}{3M\, \|
h\|_p}
\end{equation}
and hence $\| g\|_p$ should be bounded in view of condition
(\ref{eqAPRIORY}).

Thus, by Theorem 1 in \cite{LS} there is a function $g\in
L^p(\mathbb D)$ such that $g=F(g;1)$ and, consequently, by
Corollaries 2 the function $U:=N_g - {\cal P}_{N^*_g} + {\cal
P}_{\varphi}$ gives the desired solution of the Dirichlet problem
for the quasilinear Poisson equation (\ref{eqQUASILINEAR}). $\Box$

\bigskip

{\bf Corollary 3.} {\it Let $D$ be a smooth Jordan domain in
$\mathbb C$, $\Phi :\partial D\to\mathbb R$ be a continuous
function, $H: D\to\mathbb R$ be a function in the class $L^p(D)$,
$p>1$, and a function $f:\mathbb R\to\mathbb R$ be continuous and
satisfy (\ref{eqAPRIORY}).

Then there is a continuous function $u:\overline{D}\to\mathbb R$
with $u|_{\partial D}=\Phi$,
 $u\in
W^{2,p}_{\rm loc}(D)$,
\begin{equation} \label{eqQUASILINEARC}
\triangle\, u(\zeta)\ =\ H(\zeta)\cdot f(u(\zeta))\ \ \ \ \ \ \
\mbox{for a.e.}\ \zeta\in D\ .
\end{equation}
Moreover, $u\in W^{1,q}_{\rm loc}(D)$ for some $q>2$ and $u$ is
locally H\"older continuous in $D$. If in addition $\Phi$ is
H\"older continuous, then $u$ is H\"older continuous in
$\overline{D}$. Furthermore, if $p>2$, then $u\in C^{1,\alpha}_{\rm
loc}(D)$ where $\alpha = (p-2)/p$.

In particular, $u\in C^{1,\alpha}_{\rm loc}(D)$ for all
$\alpha\in(0,1)$ if $h\in L^{\infty}(D)$. If in addition $\Phi$ is
H\"older continuous on $\partial D$ with some order $\beta\in(0,1)$,
then $u$ is H\"older continuous in $\overline{ D}$ with the same
order.}

\medskip

{\bf Proof.} Let $\omega$ be a conformal mapping of $D$ onto
$\mathbb D$. By the Caratheodory-Osgood-Taylor theorem, $\omega$ is
extended to a homeomorphism $\tilde{\omega}$ of $\overline{D}$ onto
$\overline{\Bbb D}$, see \cite{C} and \cite{OT}, see also
\cite{Arsove} and Theorem 3.3.2 in the monograph \cite{CL}. Then,
setting $\varphi =\Phi\circ\tilde{\omega}^{-1}|_{\partial\mathbb
D}$, we see that the function $\varphi:\partial\mathbb D\to\mathbb
R$ is continuous. Let $h=J\cdot H\circ\Omega$ where $\Omega$ is the
inverse mapping $\omega^{-1}:\mathbb D\to D$ and $J$ is its Jacobian
$J=|\Omega^{\prime}|^2$. By the known Warschawski result, see
Theorem 2 in \cite{W}, its derivative $\Omega^{\prime}$ is extended
by continuity onto $\overline{\Bbb D}$. Consequently, $J$ is bounded
and the function $h$ is of the same class in $\mathbb D$ as $H$ in
$D$. Let $U$ be a solution of the Dirichlet problem from Proposition
1 for the equation (\ref{eqQUASILINEAR}) with the given $\varphi$
and $h$. Note that $\omega^{\prime}=1/\Omega^{\prime}\circ \omega$
is also extended by continuity onto $\overline D$ because
$\Omega^{\prime}\neq 0$ on $\partial\mathbb D$ by Theorem 1 in
\cite{W}. Thus, $u=U\circ\omega$ is the desired solution of the
Dirichlet problem for the equation (\ref{eqQUASILINEARC}). $\Box$

\section{Case of inhomogeneous and anisotropic media}

By the mentioned above  factorization theorem from \cite{GNR2017},
the study of semi--linear equations (\ref{eqQUASI}) {\it in Jordan
domains} $D$ is reduced, by means of a suitable quasiconformal
change of variables, to the study of the corresponding quasilinear
Poisson equations (\ref{QUASII}) {\it in the unit disk} $\mathbb D$.

\medskip

{\bf Theorem 4.} {\it Let $D$ be a Jordan domain in $\mathbb C$
satisfying the quasihyperbolic boun\-da\-ry condition, $A\in
M^{2\times 2}_K(D)$, $\varphi :\partial D\to\mathbb R$ be a
con\-ti\-nu\-ous function and let a function $f:\mathbb R\to\mathbb
R$ be continuous and
\begin{equation} \label{eqEQ}
\lim\limits_{t\to \infty}\ \frac{f(t)}{t}\ =\ 0\ .
\end{equation} Then there is a weak solution $u:{D}\to\mathbb
R$ of the equation (\ref{eqQUASI}) which is locally H\"older
continuous in $D$ and continuous in $\overline{D}$ with
$u|_{\partial D}=\varphi$. If in addition $\varphi$ is H\"older
continuous, then $u$ is H\"older continuous in $\overline{ D}$. }

\medskip

Following \cite{GNR2017}, under a {\bf weak solution} to the
equation (\ref{eqQUASI}) we understand a function $u\in C\cap
W^{1,2}_{\rm loc}(\Omega)$ such that, for all $\eta\in C\cap
W^{1,2}_0(D)$,
\begin{equation}\label{eqWEEK} \int\limits_D \langle A(z)\nabla
u(z),\nabla\eta(z)\rangle\ dm(z)\ +\ \int\limits_D f(u(z))\,
\eta(z)\ dm(z)\ =\ 0\ . \end{equation}

\medskip

{\bf Proof.}  By Theorem 4.1 in \cite{GNR2017}, if $u$ is a week
solution of (\ref{eqQUASI}), then $u=U\circ\omega$ where $\omega $
is a quasiconformal mapping of $D$ onto the unit disk $\mathbb D$
agreed with $A$ and $U$ is a week solution of the equation
(\ref{eqQUASILINEAR}) with $h=J$, where $J$ stands for the Jacobian
of $\omega^{-1}$. It is also easy to see that if $U$ is a week
solution of (\ref{eqQUASILINEAR}) with $h=J$, then $u=U\circ\omega$
is a week solution of (\ref{eqQUASI}). It allows us to reduce the
Dirichlet problem for equation (\ref{eqQUASI}) with continuous
boundary function $\varphi$ in the simply connected Jordan domain
$D$  to the Dirichlet problem for the equation (\ref{eqQUASILINEAR})
in the unit disk $\mathbb D$ with the continuous boundary function
$\psi=\varphi\circ\omega^{-1}$. Indeed, $\omega$ is extended to a
ho\-meo\-mor\-phism of $\overline D$ onto $\overline{\Bbb D}$, see
e.g. Theorem I.8.2 in \cite{LV:book}. Thus, the function $\psi$ is
well defined and really is continuous on the unit circle.

It is well-known that the quasiconformal mapping $\omega$ is locally
H\"older continuous in $D$, see  Theorem 3.5 in \cite{Bo^*}. Taking
into account the fact that $D$ is a Jordan domain in $\mathbb C$
satisfying the quasihyperbolic boun\-da\-ry condition, we can show
that both mappings $\omega$ and $\omega^{-1}$ are H\"older
continuous in $\overline{D}$ and $\overline{\mathbb D}$,
correspondingly. Indeed,
 $\omega =
H\circ \Omega$ where $\Omega$ is a conformal (Riemann) mapping of
$D$ onto $\mathbb D$ and $H$ is a quasiconformal mapping of $\mathbb
D$ onto itself. The mappings $\Omega$ and $\Omega^{-1}$ are H\"older
continuous in $\overline{D}$ and in $\overline{\mathbb D}$,
correspondingly, by Theorem 1 and its corollary in \cite{BP}. Next,
by the reflection principle $H$ can be extended to a quasiconformal
mapping of $\mathbb C$ onto itself, see e.g. I.8.4 in
\cite{LV:book}, and, consequently, $H$ and $H^{-1}$ are also
H\"older continuous in $\overline{\mathbb D}$, see again Theorem 3.5
in \cite{Bo^*}. The H\"older continuity of $\omega$ and
$\omega^{-1}$ in closed domains follows immediately.
\par
Now it is easy to see that if $\varphi$ is H\"older continuous, then
$\psi$ is also so, and all the conclusions of Theorem 4 follow from
Theorem 3. $\Box$

\medskip

{\bf Corollary 4.} {\it In particular, under hypotheses of Theorem 4
on $D$, $\varphi$ and $f$, there is a weak solution $U$ of the
quasilinear Poisson equation
\begin{equation} \label{eqQUASILINEARP}
\triangle\, U(z)\ =\ f(U(z))\ \ \ \ \ \ \ \mbox{for a.e.}\
z\in\mathbb D\ .
\end{equation}
which is locally H\"older continuous in $D$ and continuous in
$\overline{D}$ with $U|_{\partial D}=\varphi$. If in addition
$\varphi$ is H\"older continuous, then $U$ is H\"older continuous in
$\overline{ D}$. }

\medskip

Recall that a domain $D$ in ${\mathbb R}^n$, $n\ge 2$, is called
satisfying {\bf (A)--condition} if
\begin{equation} \label{eqA}
\hbox{mes}\ D\cap B(\zeta,\rho)\ \le\ \Theta_0\cdot\hbox{mes}\
B(\zeta,\rho)\ \ \ \ \ \ \ \ \forall\ \zeta\in\partial D\ ,\
\rho\le\rho_0
\end{equation}
for some $\Theta_0$ and $\rho_0\in(0,1)$, see 1.1.3 in \cite{LU}.
Recall also that a domain $D$ in ${\mathbb R}^n$, $n\ge 2$, is said
to be satisfying the {\bf outer cone condition} if there is a cone
that makes possible to be touched by its top to every boundary point
of $D$ from the completion of $D$ after its suitable rotations and
shifts. It is clear that the outer cone condition implies
(A)--condition.

\medskip

{\bf Remark 3.} Note that quasidisks $D$ satisfy (A)--condition.
Indeed, the quasidisks are the so--called $QED-$domains by
Gehring--Martio, see Theorem 2.22 in \cite{GM1}, and the latter
satisfy the condition
\begin{equation} \label{eqQ} \hbox{mes}\ D\cap B(\zeta,\rho)\ \ge\
\Theta_*\cdot\hbox{mes}\ B(\zeta,\rho)\ \ \ \ \ \ \ \ \forall\
\zeta\in\partial D\ ,\ \rho\le\hbox{dia} D
\end{equation}
for some $\Theta_*\in(0,1)$, see Lemma 2.13 in \cite{GM1}, and
quasidisks (as domains with quasihyperbolic boundary) have
boundaries of the Lebesgue measure zero, see e.g. Theorem 2.4 in
\cite{AK}. Thus, it remains to note that, by definition, the
completions of quasidisks $D$ in the the extended complex plane
$\overline{\mathbb C}:=\mathbb C\cup\{\infty\}$ are also quasidisks
up to the inversion with respect to a circle in $D$.

\medskip

Probably the first example of a simply connected plane domain $D$
with the quasihyperbolic boundary condition which is not a quasidisk
was constructed in \cite{BP}, Theorem 2. However, this domain had
(A)--condition. Probably one of the simplest example of a domain $D$
with the quasihyperbolic boundary condition and without
(A)--condition is the union of 3 open disks with the radius 1
centered at the points $0$ and $1\pm i$. It is clear that the domain
has zero interior angle in its boundary point $1$ and by Remark 4 it
is not quasidisk.

\bigskip

\section{The Dirichlet problem in terms of prime ends}

The simplest example of a domain $D$ with the quasihyperbolic
boundary condition and simultaneously without (A)--condition is disc
with a cut along its radius.

\medskip

Before to formulate the corresponding results for non-Jordan
domains, let us recall the necessary de\-fi\-ni\-tions of the
relevant notions and notations. Namely, we follow Caratheodory
\cite{Car$_2$} under the definition of the prime ends of domains in
$\Bbb C$, see also Chapter 9 in \cite{CL}. First of all, recall that
a continuous mapping $\sigma: \Bbb I\to \Bbb C$, $\Bbb I=(0,1)$, is
called a {\bf Jordan arc} in $\Bbb C$ if
$\sigma(t_1)\neq\sigma(t_2)$ for $t_1\neq t_2$. We also use the
notations $\sigma$, $\overline{\sigma}$ and $\partial\sigma$ for
$\sigma(\Bbb I)$, $\overline{\sigma(\Bbb I)}$ and
$\overline{\sigma(\Bbb I)}\setminus\sigma(\Bbb I)$, correspondingly.
A {\bf cross--cut} of a simply connected domain $D\subset\Bbb C$ is
a Jordan arc $\sigma$ in the domain $D$ with both ends on $\partial
D$ splitting $D$.

A sequence $\sigma_1,\ldots, \sigma_m,\ldots$ of cross-cuts of $D$
is called a {\bf chain} in $D$ if:

(i) $\overline{\sigma_i}\cap\overline{\sigma_j}=\varnothing$ for
every $i\neq j$, $i,j= 1,2,\ldots$;

(ii) $\sigma_{m}$ splits $D$ into 2 domains one of which contains
$\sigma_{m+1}$ and another one $\sigma_{m-1}$ for every $m>1$;

(iii)  $\delta(\sigma_{m})\to0$ as $m\to\infty$ where
$\delta(\sigma_{m})$ is the diameter of $\sigma_m$ with respect to
the Euclidean metric in $\mathbb C$.

Correspondingly to the definition, a chain of cross-cuts $\sigma_m$
generates a sequence of do\-mains $d_m\subset D$ such that
$d_1\supset d_2\supset\ldots\supset d_m\supset\ldots$ and $\,
D\cap\partial\, d_m=\sigma_m$. Chains of cross-cuts $\{\sigma_m\}$
and $\{\sigma_k'\}$ are called {\bf equivalent} if, for every
$m=1,2,\ldots$, the domain $d_m$ contains all domains $d_k'$ except
a finite number and, for every $k=1,2,\ldots$, the domain $d_k'$
contains all domains $d_m$ except a finite number, too. A {\bf prime
end} $P$ of the domain $D$ is an equivalence class of chains of
cross-cuts of $D$. Later on, $E_D$ denote the collection of all
prime ends of a domain $D$ and $\overline D_P=D\cup E_D$ is its
completion by its prime ends.

\medskip

Next, we say that a sequence of points $p_l\in D$ is {\bf convergent
to a prime end} $P$ of $D$ if, for a chain of cross--cuts $\{
\sigma_m\}$ in $P$, for every $m=1,2,\ldots$, the domain $d_m$
contains all points $p_l$ except their finite collection. Further,
we say that a sequence of prime ends $P_l$ converge to a prime end
$P$ if, for a chain of cross--cuts $\{ \sigma_m\}$ in $P$, for every
$m=1,2,\ldots$, the domain $d_m$ contains chains of cross--cuts $\{
\sigma_k'\}$ in all prime ends $P_l$ except their finite collection.

\medskip

A basis of neighborhoods of a prime end $P$ of $D$ can be defined in
the following way. Let $d$ be an arbitrary domain from a chain in
$P$. Denote by $d^*$ the union of $d$ and all prime ends of $D$
having some chains in $d$. Just all such $d^*$ form a basis of open
neighborhoods of the prime end $P$. The corresponding topology on
$E_D$ and, respectively, on $\overline D_P$ is called the {\bf
topology of prime ends}. The continuity of functions on $E_D$ and
$\overline D_P$ will be understood with respect to this topology or,
the same, with respect to the above convergence.

\bigskip

{\bf Theorem 5.} {\it Let $D$ be a bounded simply connected domain
in $\mathbb C$ satisfying the quasihyperbolic boun\-da\-ry
condition, $A\in M^{2\times 2}_K(D)$, $\varphi :E_D\to\mathbb R$ be
a con\-ti\-nu\-ous function and let a function $f:\mathbb
R\to\mathbb R$ be continuous and
\begin{equation} \label{eqEQE}
\lim\limits_{t\to +\infty}\ \frac{f(t)}{t}\ =\ 0\ .
\end{equation} Then there is a weak solution $u:{D}\to\mathbb
R$ of the equation (\ref{eqQUASI}) which is locally H\"older
continuous in $D$ and continuous in $\overline{D}_P$ with
$u|_{E_D}=\varphi$.}

\bigskip

{\bf Proof.}  Again by Theorem 4.1 in \cite{GNR2017}, if $u$ is a
week solution of (\ref{eqQUASI}), then $u=U\circ\omega$ where
$\omega $ is a quasiconformal map of $D$ onto the unit disk $\mathbb
D$ agreed with $A$ and $U$ is a week solution of the equation
(\ref{eqQUASILINEAR}) with $h=J$, the Jacobian of $\omega^{-1}$.
Similarly, if $U$ is a week solution of (\ref{eqQUASILINEAR}) with
$h=J$, then  $u=U\circ\omega$ is a week solution of (\ref{eqQUASI}).

Hence the Dirichlet problem for (\ref{eqQUASI}) in the domain $D$
will be reduced to the so for (\ref{eqQUASILINEAR}) in $\mathbb D$
with the corresponding boundary function
$\psi=\varphi\circ\omega^{-1}$. The existence and continuity of the
boundary function $\psi$ in the case of an arbitrary bounded simply
connected domain $D$ is a fundamental result of the theory of the
boundary behavior of conformal and quasiconformal mappings. Namely,
$$\omega^{-1}=H\circ \Omega$$
where $\Omega$ stands for a quasiconformal automorphism of the unit
disk ${\Bbb D}$ and $H$ is a conformal mapping of $\Bbb D$ onto
$\Omega.$ It is known that $\Omega$ can be extended to a
homeomorphism of $\overline{\Bbb D}$ onto itself, see e.g. Theorem
I.8.2 in \cite{LV:book}. Moreover, by the well-known Caratheodory
theorem on the boundary correspondence under conformal mappings, see
e.g. Theorems 9.4 and 9.6 in \cite{CL}, the mapping $H$ is extended
to a homeomorphism of $\overline{\Bbb D}$ onto $\overline{D}_P$.
Thus, the function $\psi$ is well defined and really continuous on
the unit circle.

Moreover, $\omega$ is locally H\"older continuous in $\mathbb D$,
see e.g. Theorem 3.5 in \cite{Bo^*}. Thus,  by Theorem 4.1 in
\cite{GNR2017}, Theorem 5 follows from Theorem 3. $\Box$

\medskip

{\bf Corollary 5.} {\it In particular, under hypotheses of Theorem 5
on $D$, $\varphi$ and $f$, there is a weak solution $U$ of the
quasilinear Poisson equation (\ref{eqQUASILINEARP}) which is locally
H\"older continuous in $D$ and continuous in  $\overline{D}_P$ with
$U|_{E_D}=\varphi$. }

\medskip

By the way, it is turned out to be possible to solve the Dirichlet
problem even for the degenerate Beltrami equations in terms of
Caratheodory prime ends in arbitrary bounded finitely connected
domains $D$ in $\mathbb C$, see \cite{GRY}.

\medskip

\section{Some applied corollaries}

The interest in this subject is well known both from a purely
theoretical point of view, due to its deep relations to linear and
nonlinear harmonic analysis, and because of numerous applications of
equations of this type  in various areas of physics, differential
geometry, logistic problems etc., see e.g. \cite{BK}, \cite{Diaz},
\cite{GT}, \cite{HKM}, \cite{KP}, \cite{Landis}, \cite{MV},
\cite{NN}, \cite{Ver} and the references therein. In particular, in
the excellent book by M. Marcus and L. Veron \cite{MV} the reader
can find a comprehensive analysis of the Dirichlet problem for the
semi--linear equation
\begin{equation}\label{Poisson3}
\triangle u(z)\ =\ f(z,u(z))
\end{equation}
in smooth ($C^2$) domains $D$ in $\mathbb R^n$, $n\ge 3$, with
boundary data in $L^1.$ Here $t \to f(\cdot , t) $ is a continuous
mapping from $\mathbb R$ to a weighted Lebesgue's space
$L^1(D,\rho)$, and $z \to f(z,\cdot)$ is a non-decreasing function
for every $z\in D$, $f(z,0)\equiv 0$, such that
\begin{equation}\label{LIM}
\lim_{t\to \infty}{f(z,t)\over t}\ =\ \infty
\end{equation}
uniformly with respect to the parameter $z$ in compact subsets of
$D$.

The mathematical modelling of some reaction-diffusion problems leads
to the study of the corresponding Dirichlet problem for the equation
(\ref{Poisson1}) with specified right hand side. Following
\cite{Aris},  a nonlinear system can be obtained for the density $u$
and the temperature $T$ of the reactant. Upon eliminating $T$ the
system can be reduced to a scalar problem for the concentration
\begin{equation}\label{RDP}
 \triangle  u\ =\ \lambda\cdot f(u)
 \end{equation}
where $\lambda$ stands for a positive constant.

It turns out that the density of the reactant $u$ may be zero in a
closed interior region $D_0$ called {\it a dead core.} If, for
instance in the equation (\ref{RDP}), $f(u) = u^q,$  $q > 0,$ a
particularization of the results in Chapter 1 of \cite{Diaz} shows
that a dead core may exist if and only if $0 < q < 1$ and $\lambda$
is large enough. See also the corresponding examples of dead cores
in \cite{GNR2017}. We have by Theorem 4 the following:

\medskip

{\bf Theorem 6.} {\it Let $D$ be a Jordan domain in $\mathbb C$
satisfying the quasihyperbolic boun\-da\-ry condition, $A\in
M^{2\times 2}_K(D)$, $\varphi :\partial D\to\mathbb R$ be a
con\-ti\-nu\-ous function. Then there is a weak solution
$u:{D}\to\mathbb R$ of the semilinear equation
\begin{equation} \label{eqQUASIQ} {\rm div\,}[\, A(z)\, \nabla u(z)\, ]\ =\
u^q(z)\ , \ \ \ 0\ <\ q\ <\ 1\ ,
\end{equation}
which is locally H\"older continuous in $D$ and continuous in
$\overline{D}$ with $u|_{\partial D}=\varphi$. If in addition
$\varphi$ is H\"older continuous, then $u$ is H\"older continuous in
$\overline{ D}$. }

\medskip

Recall that under a weak solution to the equation (\ref{eqQUASIQ})
we understand a function $u\in C\cap W^{1,2}_{\rm loc}(\Omega)$ such
that, for all $\eta\in C\cap W^{1,2}_0(D)$,
\begin{equation}\label{eqDIFFUSION} \int\limits_D \langle\, A(z)\nabla
u(z),\nabla\eta(z)\, \rangle\ dm(z)\ +\ \int\limits_D u^q(z)\,
\eta(z)\ dm(z)\ =\ 0\ . \end{equation}

We have also the following significant consequence of Corollary 3.

\medskip

{\bf Corollary 6.} {\it Let $D$ be a smooth Jordan domain in
$\mathbb C$ and $\varphi :\partial D\to\mathbb R$ be a
con\-ti\-nu\-ous function. Then there is a weak solution $U$ of the
quasilinear Poisson equation
\begin{equation}\label{eqQUASILINEARQ}
\triangle\, U(z)\ =\ U^q(z)\ , \ \ \ 0\ <\ q\ <\ 1\ ,
\end{equation}
which is continuous in $\overline{D}$ with $U|_{\partial D}=\varphi$
and $U\in C^{1,\alpha}_{\rm loc}(D)$ for all $\alpha\in(0,1)$. If in
addition $\varphi$ is H\"older continuous with some order
$\beta\in(0,1)$, then $U$ is also H\"older continuous in $\overline{
D}$ with the same order.}

\medskip

Recall also that certain mathematical models of a heated plasma lead
to nonlinear equations of the type (\ref{RDP}). Indeed, it is known
that some of them have the form $\triangle\psi(u)=f(u)$ with
$\psi'(0)=+\infty$ and $\psi'(u)>0$ if $u\not=0$ as, for instance,
$\psi(u)=|u|^{q-1}u$ under $0 < q < 1$, see e.g. \cite{Bear},
\cite{Be} and \cite{Diaz}, p. 4. With the replacement of the
function $U=\psi(u)=|u|^q\cdot {\rm sign}\, u$, we have that $u =
|U|^Q\cdot {\rm sign}\, U$, $Q=1/q$, and, with the choice $f(u) =
|u|^{q^2}\cdot {\rm sign}\, u$, we come to the equation $\triangle U
= |U|^q\cdot {\rm sign}\, U=\psi(U)$.

Of course, the similar results can be formulated for this case, for
instance:

\medskip

{\bf Corollary 7.} {\it Let $D$ be a smooth Jordan domain in
$\mathbb C$ and $\varphi :\partial D\to\mathbb R$ be a
con\-ti\-nu\-ous function. Then there is a weak solution $U$ of the
quasilinear Poisson equation
\begin{equation}\label{eqPLASMA}
\triangle\, U(z)\ =\ |U(z)|^{q-1}U(z)\ , \ \ \ 0\ <\ q\ <\ 1\ ,
\end{equation}
which is continuous in $\overline{D}$ with $U|_{\partial D}=\varphi$
and $U\in C^{1,\alpha}_{\rm loc}(D)$ for all $\alpha\in(0,1)$. If in
addition $\varphi$ is H\"older continuous with some order
$\beta\in(0,1)$, then $U$ is also H\"older continuous in $\overline{
D}$ with the same order. }

\medskip

In the combustion theory, see e.g. \cite{Barenblat},
\cite{Pokhozhaev} and the references therein, the following model
equation
\begin{equation}
{\partial u(z,t)\over \partial t}\ =\ {1\over \delta}\cdot \triangle
u\ +\ e^{u}\ ,\ \ \ t\geq 0,\ z\in D,
\end{equation}
occupies a special place. Here $u\ge 0$ is the temperature of the
medium and $\delta$ is a certain positive parameter. We restrict
ourselves by stationary solutions of the equation and its
ge\-ne\-ra\-li\-za\-tions in anisotropic and inhomogeneous media
although our approach makes it possible to consider the parabolic
case, see \cite{GNR2017}. Namely, by Theorem 4 we have the following
statement:


{\bf Theorem 7.} {\it Let $D$ be a Jordan domain in $\mathbb C$
satisfying the quasihyperbolic boun\-da\-ry condition, $A\in
M^{2\times 2}_K(D)$, $\varphi :\partial D\to\mathbb R$ be a
con\-ti\-nu\-ous function. Then there is a weak solution
$U:{D}\to\mathbb R$ of the semilinear equation
\begin{equation} \label{eqQUASIQA} {\rm div\,}[\, A(z)\, \nabla U(z)\, ]\ =\
\delta\cdot e^{-U(z)}
\end{equation}
which is locally H\"older continuous in $D$ and continuous in
$\overline{D}$ with $u|_{\partial D}=\varphi$. If in addition
$\varphi$ is H\"older continuous, then $u$ is H\"older continuous in
$\overline{ D}$. }


Finally, we obtain also the following consequence of Corollary 3.


{\bf Corollary 8.} {\it Let $D$ be a smooth Jordan domain in
$\mathbb C$ and $\varphi :\partial D\to\mathbb R$ be a
con\-ti\-nu\-ous function. Then there is a weak solution $U$ of the
quasilinear Poisson equation
\begin{equation}\label{eqQUASILINEARQ}
\triangle\, U(z)\ =\ \delta\cdot e^{U(z)}
\end{equation}
which is continuous in $\overline{D}$ with $U|_{\partial D}=\varphi$
and $U\in C^{1,\alpha}_{\rm loc}(D)$ for all $\alpha\in(0,1)$. If in
addition $\varphi$ is H\"older continuous with some order
$\beta\in(0,1)$, then $U$ is also H\"older continuous in $\overline{
D}$ with the same order.}

\medskip

Thus,  the results on regular solutions for the quasilinear Poisson
equations (\ref{eqQUASILINEAR}) and the comprehensively developed
theo\-ry of quasiconformal mappings in the plane, see e.g. the
monographs \cite{Ahlfors:book}, \cite{BGMR}, \cite{GRSY},
\cite{LV:book} and \cite{MRSY}, are  good basiс tools for the
further study of equations (\ref{eqQUASI}). The latter opens up a
new approach to the study of a number of semi-linear equations of
mathematical physics in anisotropic and inhomogeneous media.

\bigskip

\noindent {\it Institute of Applied Mathematics and Mechanics of
National Academy\\ of Sciences of Ukraine, 84100, Ukraine,
Slavyansk, 1 Dobrovolskogo Str.}

\noindent Email: vgutlyanskii@gmail.com, star-o@ukr.net,
Ryazanov@nas.gov.ua

\bigskip

\noindent{\it National University of Cherkasy, Physics Department,\\
Physics Department, Laboratory of Mathematical Physics, \\
18001, Ukraine, Cherkasy, 81 Shevchenko Blvd.}

\noindent Email: vl.ryazanov1@gmail.com


\vskip 2mm
\begin{thebibliography}{100}
\small

\bibitem{Ahlfors:book}
L.V. Ahlfors, {\em Lectures on quasiconformal mappings}, Van
Nostrand Ma\-the\-ma\-ti\-cal Studies {\bf 10}, Van Nostrand Co.,
Inc., Toronto, Ont.-New York-London, 1966; transl. as {\em Lektsii
po kvazikonformnym otobrazheniyam}, Mir, Moscow, 1969.

\bibitem{AIM-book}
K. Astala, T. Iwaniec, G. Martin, {\em Elliptic partial differential
equations and quasiconformal mappings in the plane}, Princeton
Mathematical Series {\bf 48}, Princeton University Press, Princeton,
NJ, 2009.

\bibitem{AK}
K. Astala, P. Koskela, \emph{Quasiconformal mappings and global
integrability of the derivative}, J. Anal. Math. {\bf 57} (1991),
203--220.

\bibitem{Aris}
R. Aris, {\em The Mathematical Theory of Diffusion and Reaction in
Permeable Ca\-ta\-lysts}, Volumes I and II, Clarendon Press, Oxford,
1975.

\bibitem{Arsove}
M.G. Arsove, \emph{The Osgood-Taylor-Caratheodory theorem}, Proc.
Amer. Math. Soc. {\bf 19} (1968), 38--44.

\bibitem{Barenblat}
G.I. Barenblatt, Ja.B. Zel'dovic, V.B. Librovich, G.M. Mahviladze,
\emph{Matematicheskaya teoriya goreniya i vzryva,} Nauka, Moscow,
1980; transl. in \emph{The mathematical theory of combustion and
explosions}, Consult. Bureau, New York, 1985.

\bibitem{Bear}
J. Bear, \emph{Modeling phenomena of flow and transport in porous
media}, Theory and Applications of Transport in Porous Media {\bf
31}, Springer, Cham, 2018.

\bibitem{Be}
J. Bear, \emph{Dynamics of Fluids in Porous Media}, Elsevier, New
York, 1972.

\bibitem{BP}
J. Becker, Ch. Pommerenke, \emph{H\"older continuity of conformal
mappings and nonquasiconformal Jordan curves}, Comment. Math. Helv.
{\bf 57} (1982), no. 2, 221--225.


\bibitem{Bo}
B.V. Bojarski, \emph{Homeomorphic solutions of Beltrami systems}
(Russian), Dokl. Akad. Nauk SSSR (N.S.) {\bf 102} (1955), 661–664.

\bibitem{Bo^*}
B.V. Bojarski, \emph{Generalized solutions of a system of
differential equations of the first order and elliptic type with
discontinuous coefficients}, Report Dept. Math. Stat. {\bf 118},
Univ. of Jyv\"askyl\"a, Jyv\"askyl\"a (2009); transl. from Mat. Sb.
(N.S.) {\bf 43(85)} (1957), 451–503.

\bibitem{BGMR}
B. Bojarski, V. Gutlyanskii, O. Martio, V. Ryazanov, {\em
Infinitesimal geometry of quasiconformal and bi-lipschitz mappings
in the plane}, EMS Tracts in Ma\-the\-ma\-tics {\bf 19}, European
Mathematical Society, Z\"urich, 2013.

\bibitem{BI}
{\it B. Bojarski, T. Iwaniec,} Analytical foundations of the theory
of quasiconformal mappings in Rn // Ann. Acad. Sci. Fenn. Ser. A I
Math. 8 (1983), no. 2, 257–324.

\bibitem{BK}
M. Borsuk, V. Kondratiev,  {\em Elliptic boundary value problems of
second order in piecewise smooth domains}, North-Holland
Mathematical Library {\bf 69}, Elsevier Science, Amsterdam, 2006.

\bibitem{C}
C. Caratheodory, \emph{Uber die gegenseitige Beziehung der Rander
bei der konformen Abbildungen des Inneren einer Jordanschen Kurve
auf einen Kreis}, Math. Ann. {\bf 73} (1913), 305-320.

\bibitem{Car$_2$} C. Caratheodory {\em \"Uber die Begrenzung der einfachzusammenh\"angender
Gebiete} // Math. Ann. \textbf{73} (1913),  323--370.

\bibitem{CL} E.F. Collingwood, A.J. Lohwator, {\em The Theory of Cluster
Sets}, Cambridge Tracts in Math. and Math. Physics \textbf{56},
Cambridge Univ. Press, Cambridge, 1966.

\bibitem{Diaz}
J.I. Diaz, {\em Nonlinear partial differential equations and free
boundaries}, Volume I Elliptic equations, Research Notes in
Mathematics {\bf 106}, Pitman, Boston, 1985.

\bibitem{DS}
N. Dunford, J.T. Schwartz, {\em Linear Operators. I. General
Theory}, Pure and Applied Mathematics {\bf 7}, Interscience
Publishers, New York, London, 1958; transl. as {\it Lineinye
operatory. Chast' I: Obshchaya teoriya}, Inostran. Lit., Moscow,
1962.

\bibitem{Fe}
{H. Federer}, {\em Geometric Measure Theory}, Springer-Verlag,
Berlin, 1969; transl. as {\it Geometricheskaya teoriya mery}, Nauka,
Moscow, 1987.

\bibitem{Ge}
F.W. Gehring, \emph{Characteristic properties of quasidisks},
Seminaire de Mathґematiques Superieures {\bf 84} Presses de
l’Universitґe de Montrґeal, Montreal, Que., 1982.

\bibitem{GH}
F.W. Gehring, K. Hag, {\em The ubiquitous quasidisk}, Mathematical
Surveys and Monographs {\bf 184}, American Mathematical Society,
Providence, RI, 2012.

\bibitem{GM1}
F.W. Gehring, O. Martio, \emph{Quasiextremal distance domains and
extension of quasiconformal mappings}, J. Analyse Math. {\bf 45}
(1985), 181–206.

\bibitem{GM2}
{\it Gehring F.W., Martio O.} Lipschitz classes and quasiconformal
mappings // Ann. Acad. Sci. Fenn. Ser. A I Math. -- 1985. -- {\bf
10}. -- P. 203–219.

\bibitem{GM$^*$}
M. Giaquinta, L. Martinazzi, {\em An introduction to the regularity
theory for elliptic systems, harmonic maps and minimal graphs},
Second edition. Appunti. Scuola Normale Superiore di Pisa (Nuova
Serie) [Lecture Notes. Scuola Normale Superiore di Pisa (New
Series)], 11. Edizioni della Normale, Pisa, 2012.

\bibitem{GT}
D. Gilbarg, N. Trudinger, {\em Elliptic partial differential
equations of second order}, Grundlehren der Mathematischen
Wissenschaften {\bf 224}, Springer-Verlag, Berlin, 1983; transl. as
{\em Ellipticheskie differentsial'nye uravneniya s chastnymi
proizvodnymi vtorogo poryadka}, Nauka, Moscow, 1989.


\bibitem{GNR2017}
V.Ya. Gutlyanskii, O.V. Nesmelova, V.I. Ryazanov, \emph{ On
quasiconformal maps and semi-linear equations in the plane}, Ukr.
Mat. Visn. {\bf 14} (2017), no. 2, 161-191; transl. in J. Math. Sci.
(US) {\bf 229} (2018), no. 1, 7-29.

\bibitem{GRSY}
V. Gutlyanskii, V. Ryazanov, U. Srebro, E. Yakubov, {\em The
Beltrami Equation: A Geometric Approach}, Developments in
Mathematics {\bf 26}, New York etc., Springer, 2012.

\bibitem{GRY}
V.Ya. Gutlyanskii, V.I. Ryazanov, E. Yakubov, \emph{The Beltrami
equations and prime
 ends},  Ukr.
Mat. Visn. {\bf 12} (2015), no. 1, 27--66; transl. in J. Math. Sci.
(N.Y.) {\bf 210} (2015), no. 1, 22–-51.

\bibitem{HKM}
J. Heinonen, T. Kilpelдinen, O. Martio, {\em Nonlinear potential
theory of degenerate elliptic equations}, Clarendon Press, Oxford,
New York, Tokyo, 1993.

\bibitem{Ho}
L. H\"ormander, {\em Notions of convexity}, Progress in Mathematics
{\bf 127}, Birkh\"auser Boston Inc, Boston, 1994.

\bibitem{Hor}
L. H\"ormander, {\em The analysis of linear partial differential
operators. V. {\rm I}. Distribution theory and Fourier analysis},
Grundlehren der Mathematischen Wissenschaften {\bf 256},
Springer-Verlag, Berlin, 1983; transl. as {\em Analiz lineinykh
differentsial'nykh operatorov s chastnymi proizvodnymi}. T. 1, Mir,
Moscow, 1986.

\bibitem{Ka}
D. Kalaj, \emph{ On some integral operators related to the Poisson
equation}, Integral Equations Operator Theory {\bf 72} (2012), no.
4, 563–575.

\bibitem{KV}
D. Kalaj, D. Vujadinoviґc, \emph{ The solution operator of the
inhomogeneous Dirichlet problem in the unit ball}, Proc. Amer. Math.
Soc. {\bf 144} (2016), no. 2, 623–635.

\bibitem{Ko} {P. Koosis}, {\em Introduction to $H^p$ spaces}, {Cambridge Tracts in
Ma\-the\-ma\-tics} \textbf{115}, Cambridge Univ. Press, Cambridge,
1998.

\bibitem{KP}
I. Kuzin, S. Pohozaev, {\em Entire solutions of semilinear elliptic
equations}, Progress in Nonlinear Differential Equations and their
Applications \textbf{33}, Birkh\"auser Verlag, Basel, 1997.

\bibitem{Landis}
E.M. Landis, {\em Second order equations of elliptic and parabolic
type}, Translations of Mathematical Monographs \textbf{171},
American Mathematical Society, Providence, RI, 1998.

\bibitem{La}
N.S. Landkof, {\em Foundations of modern potential theory}, Nauka,
Moscow, 1966; transl. in Grundlehren der mathematischen
Wissenschaften {\bf 180}, Springer-Verlag, New York--Heidelberg,
1972.

\bibitem{LU}
O.A. Ladyzhenskaya, N.N. Ural'tseva, {\em Linear and quasilinear
elliptic equations}, Academic Press, New York-London, 1968; transl.
from {\em Lineinye i kvazilineinye uravneniya иllipticheskogo tipa},
Nauka, Moscow 1964.

\bibitem{LV:book}
O. Lehto, K.I. Virtanen, {\em Quasiconformal mappings in the plane,}
2-nd Edition, Springer--Verlag, Berlin - Heidelberg -  New York,
1973.

\bibitem{LS}
J. Leray, Ju. Schauder, \emph{ Topologie et equations
fonctionnelles} (French), Ann. Sci. Ecole Norm. Sup. {\bf 51}
(1934), no. 3, 45--78; transl. as \emph{ Topology and functional
equations}, Uspehi Matem. Nauk (N.S.) {\bf 1} (1946), no. 3-4
(13-14), 71--95.


\bibitem{MV}
M. Marcus, L. Veron, {\em Nonlinear second order elliptic equations
involving measures}, De Gruyter Series in Nonlinear Analysis and
Applications {\bf 21}, De Gruyter, Berlin, 2014.

\bibitem{MRSY}
O. Martio, V. Ryazanov, U. Srebro, E. Yakubov, {\em Moduli in Modern
Mapping Theory}, New York, Springer, 2009.

\bibitem{NN}
R.M. Nasyrov, S.R. Nasyrov,   \emph{ Convergence of S. A.
Khristianovich's approximate method for solving the Dirichlet
problem for an elliptic equation} (Russian) Dokl. Akad. Nauk SSSR
{\bf 291} (1986), no. 2, 294-298.

\bibitem{OT}
W. Osgood, E. Taylor, \emph{Conformal transformations on the
boundaries of their regions of definition}, Trans. Amer. Math. Soc.
{\bf 14} (1913), 277-298.

\bibitem{Pokhozhaev}
S.I. Pokhozhaev, \emph{ On an equation of combustion theory}, Mat.
Zametki {\bf 88} (2010), no. 1, 53--62; transl. in Math. Notes {\bf
88} (2010), no. 1-2, 48–56.


\bibitem{Ra}
T. Ransford, {\em Potential theory in the complex plane}, London
Mathematical Society Student Texts {\bf 28}, Cambridge University
Press, Cambridge, 1995.

\bibitem{Ru}
W. Rudin, {\em Real and complex analysis}, McGraw-Hill Book Co., New
York, 1987.

\bibitem{Sa} {S. Saks}, {\em Theory of the integral}, Warsaw,
1937; Dover Publications Inc., New York, 1964; transl. as {\em
Teoria integrala}, Inostran. Lit., Moscow, 1949.


\bibitem{So}
S.L. Sobolev, {\em Some applications of functional analysis in
mathematical physics}, Transl. Math. Mon. {\bf 90}, AMS, Providence,
RI, 1991; from {\em Nekotorye primeneniya funkcional'nogo analiza v
matematiceskoi fizike}, Nauka, Moscow, 1988.

\bibitem{Ve} { I.N. Vekua,} \emph{ Generalized analytic functions},
Pergamon Press, London-Paris-Frankfurt; Addison-Wesley Publishing
Co., Inc., Reading, Mass. 1962; transl. from {\em Obobshchennye
analiticheskie funktsii}, Fiz.-Mat. Lit., Moscow 1959.

\bibitem{Ver}
L. Veron, \emph{ Local and global aspects of quasilinear degenerate
elliptic equations. Quasilinear elliptic singular problems}, World
Scientific Publishing Co. Pte. Ltd., Hackensack, NJ, 2017.

\bibitem{W}
S.E. Warschawski, {\em On differentiability at the boundary in
conformal mapping}, Proc. Amer. Math. Soc. {\bf 12} (1961), 614–620.




\end{thebibliography}
\end{document}